\documentclass[preprint,12pt]{elsarticle}

\usepackage{xcolor, graphicx}
\usepackage{amssymb}
\usepackage{ifpdf}
\usepackage{diagbox}
\usepackage{amsmath}
\usepackage{verbatim}

\makeatletter
\let\elsarticle@keyword\keyword
\g@addto@macro\frontmatter{\let\keyword\elsarticle@keyword}
\makeatother
\usepackage{program}

\journal{Discrete Mathematics}

\begin{document}

\begin{frontmatter}



\title{Equivalence Classes Induced by Binary Tree Isomorphism -- Generating Functions}


\author[x1]{David Serena}
\address[x1]{CTO, CogNueva, Inc., California, USA\\ info@cognueva.com}
%\ead[x1]{info@cognueva.com}

\author[x2]{William J Buchanan}
\address[x2]{Blockpass ID Lab, Edinburgh Napier University\\ b.buchanan@napier.ac.uk}
%\ead[x2]{b.buchanan@napier.ac.uk}

\begin{abstract}
Working with generating functions, the combinatorics of a recurrence relation can be expressed in a way that allows
for more efficient calculation of the quantity.  This is true of the Catalan Numbers for an ordered binary tree.   
Binary tree isomorphism is an important problem in computer science.  
The enumeration of the number of non-isomorphic rooted binary trees is therefore well known.  The paper reiterates the 
known results for ordered binary trees and presents previous results for enumeration of non-isomorphic rooted binary trees.   Then new enumeration results are put forward for the two-color binary tree isomorphism parametrized by the number of nodes, the number of specific color and the number of non-isomorphic sibling subtrees.   Multi-variate generating function equations are presented that enumerate these tree structures.  The generating functions with these parameterizations separate multiplicatively into simplified generating function equations. 
\end{abstract}

\begin{keyword}
binary tree isomorphism \sep binary tree \sep graph isomorphism \sep combinatorics
\end{keyword}

\end{frontmatter}

% \linenumbers


\section{Introduction}

In Sections \ref{OBT} and  \ref{NIBT} the basic known results are reiterated with derivation.  Then in Sections \ref{NNIRBTWNNISG}, \ref{TCRNIBT} and \ref{TCRBTIPNSCNISN} new enumeration results are put forward for two-color binary tree isomorphism parametrized by number of nodes \cite{jaffke2024b}, number of specific color and number of non-isomorphic sibling subtrees.   Multi-variate generating function equations are presented that enumerate these tree structures.  The generating functions with these parameterizations separate multiplicatively into simplified generating function equations. 

\section{Basic Ordered Binary Trees}
\label{OBT}

Starting with the most basic result.  When binary trees are enumerated with regard to the ordering of the subtrees, the standard enumeration of ordered binary trees\footnote{Ordered trees are where the 
left and right position of the nodes matters.}  is as follows:

Such trees are given by the recurrence when $n\geq 1$.

\[
C_n=\sum_{k=0}^{n-1}  C_{n-1-k } C_k
\]

The base case clearly shows that $C_0=1$. When the following generating function
$R(z)$ is defined as

\[
F(z) = \sum_{n=0}^\infty C_n z^n
\]

Substituting the recurrence yields the following equation.

\begin{eqnarray}
\label{GFE0} z F(z)^2= F(z) -1
\end{eqnarray}

Since the equation is a quadratic it yields a closed form solution for $R(z)$. 

\[
F(z)={{1-{\sqrt{1-4z}}}\over {2z}}
\]

Which allows for a direct closed form expression of the coefficient in terms
of the binomial function.  These are just the Catalan Numbers.

\[
C_n= {1 \over {n+1}} \biggl( {{2n}\atop{n}} \biggr)
\]

\section{Non-Isomorphic Binary Trees With Labeled Root Node}
\label{NIBT}

The following is a rederivation of the recurrence which represents the number of non-isomorphic binary trees with labeled root nodes and $n$ nodes in the
tree. \cite{Etherington,etherington1937non,etherington1940some}  An alternate way of looking at the 
isomorphism problem with binary trees is the concept of ``flip-equivalence'' \cite{bergold2023topological}.

Can we map one binary tree with a designated root to another by flipping the children of each node?    This is equivalent to the isomorphism of the trees where the root node is labeled.

Building the recurrence, the base cases are trivially. 

\[
B_0=B_1=1
\]

The recursive cases are distinct for even and odd $n$ respectively ($n\geq 1$),  

\begin{eqnarray*}
B_{2 n} &=&  {1 \over 2} \sum_{k=0}^{2n-1} B_{2n-k-1} B_k \\
B_{2 n+1} &=& {1 \over 2}\left(\left[ \sum_{k=0}^{2n} B_{2n-k} B_k \right] +
  B_n \right) 
\end{eqnarray*}

\begin{figure}[h]
\begin{center}
\resizebox{!}{4.5in}{
\ifpdf
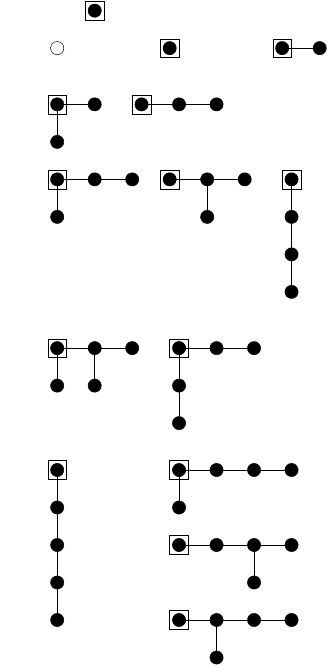
\else
\input{Etherington.pstex_t}
\fi
}
\end{center}
\caption{Rooted Non-Isomorphic Binary Trees}
\label{rbt}
\end{figure}

The first equation is derived by observing that when there are an even number
of nodes in a tree, the subtrees, possibly empty, formed by the children of the
root may never be isomorphic by a simple counting argument. As order does not
matter with regards to isomorphism, the factor of ${1\over 2}$ appears outside
the summation.

The second equation is more complex to interpret.  If the number of nodes in a
tree is odd there is a possibility that the children of the root contain the
same number of nodes.  Outside of the given base cases, there are two
possibilities.  Either the siblings are the roots of isomorphic subtrees or
they are not. The former case is handled by the first term and the latter is
handled by the second term. As any tree with more than two nodes has at least
two non-isomorphic manifestations, all the cases are covered.   Figure \ref{rbt} shows the 
enumeration of these trees.  

\subsection{Generating Functions}

\begin{eqnarray*}
G(s) &=& \sum_{k=0}^\infty B_k s^k\\
H(s) &=& \sum_{k=0}^\infty  B_{2k} s^k\\
I(s) &=& \sum_{k=0}^\infty  B_{2k+1} s^k 
\end{eqnarray*}

\subsection{Equation in G(s)}
\label{SIMPLEGEN}

\begin{eqnarray*}
G(s) &=& 1 + s + \sum_{n=2}^\infty B_n s^n \\
&=& 1+s + \sum_{n=1}^\infty (B_{2n} s^{2n} + B_{2n+1} s^{2n+1})\\
&=& 1+s + {1 \over 2} \sum_{n=1}^\infty \biggl(\sum_{\ell=0}^{2n-1}
B_{2n-\ell-1} 
B_\ell s^{2n} + \left[ \sum_{k=0}^{2n} B_{2n-k} B_k s^{2n+1} \right] +
B_n s^{2n+1}\biggr)\\
&=&  1+{s\over 2} + {1 \over 2} \biggr[\sum_{n=1}^\infty \biggl[\sum_{\ell=0}^{n-1} ( \chi(2
\vert n) B_{n-\ell-1}
B_\ell +  \chi(2 \nmid n)B_{n-\ell-1} B_\ell)\biggr]s^{n}  +
B_n s^{2n+1}\biggl]
\end{eqnarray*}

Where $\chi(p)=1$ if $p$ is true or $\chi(p)=0$ if $p$ is false.    $k \vert n$ is true if $k$ evenly divides $n$ and false otherwise.  $k \nmid n$ is the negation $k \vert n$. 

\[
G(s)= 1 +{s\over 2}+ {1\over 2} \sum_{n=1}^\infty \sum_{k=0}^{n-1} B_{n-k-1} B_k s^n
+  {1 \over 2} \sum_{n=1}^\infty B_n s^{2n+1}
\]

Changing the indices

\[
= 1 +  {1\over 2} \sum_{n=0}^\infty  \sum_{k=0}^\infty B_n B_k
s^{n+k+1} + {s \over 2} \sum_{n=0}^\infty B_n s^{2n}
\]

Since $ \sum_{n=0}^\infty  \sum_{k=0}^\infty B_n B_k
s^{n+k+1}= s \left(\sum_{n=0}^\infty  B_n s^n \right) \left( \sum_{k=0}^\infty
B_k s^{k}\right)=s(G(s))^2$ then the generating function is given by 

\begin{eqnarray*}
G(s) = 1 + {s \over 2}\biggl[G(s)^2 + G(s^2)\biggr]
\end{eqnarray*}

\begin{eqnarray*}
G(s)&=&1 + s + {s^2} + 2\,{s^3} + 3\,{s^4} + 6\,{s^5} + 11\,{s^6}  +
  23\,{s^7} + 46\,{s^8} + 98\,{s^9} + 207\,{s^{10}} \\ 
&& + 451\,{s^{11}} + 983\,{s^{12}} + 2179\,{s^{13}} + 
  4850\,{s^{14}} + 10905\,{s^{15}}+\cdots
\end{eqnarray*}

\begin{eqnarray*}
G(s)^2&=&1 + 2\,s + 3\,{s^2} + 6\,{s^3} + 11\,{s^4} + 22\,{s^5} + 
  44\,{s^6} + 92\,{s^7} + 193\,{s^8} + 414\,{s^9} \\&& + 
  896\,{s^{10}}+ 1966\,{s^{11}} + 4347\,{s^{12}} + 
  9700\,{s^{13}} + 21787\,{s^{14}} + 49262\,{s^{15}} + \cdots
\end{eqnarray*}

\[
G(s)=H(s^2)+s I(s^2)
\]

Similarly, equations for $H(s)$ and $I(s)$ can be developed.

\begin{eqnarray*}
H(s)&=& 1+ {1\over 2} \sum_{n=1}^\infty \sum_{k=0}^{2n-1} B_{2n-k-1} B_k s^n\\
I(s)&=& 1+ {1\over 2} \sum_{n=1}^\infty \sum_{k=1}^{2n}  B_{2n-k} B_k s^n + {1
  \over 2} \sum_{n=1}^\infty B_n s^n
\end{eqnarray*}

\begin{eqnarray*}
H(s)&=& 1 + {1\over 2} \sum_{n=1}^\infty \left( \sum_{k=0}^{n-1} B_{2n-1-2k}
  B_{2k} s^n + B_{2(n-k-1)} B_{2k+1} \right) s^n \\
&=& 1 + {s\over 2} \sum_{n=0}^\infty \sum_{k=0}^\infty \left(B_{2n+1} B_{2k} 
  + B_{2n} B_{2k+1}\right ) s^{n+k}\\
&=& 1+ {s\over 2}  \biggl[ \left(\sum_{n=0}^\infty B_{2n+1} s^n \right)
  \left( \sum_{k=0}^\infty B_{2k} s^k\right) 
  + \left(\sum_{n=0}^\infty B_{2n} s^n \right) \left( \sum_{k=0}^\infty
  B_{2k+1} 
  s^k \right) \biggr]\\
&=& 1+ s I(s) H(s) 
\end{eqnarray*}

\begin{eqnarray*}
I(s) &=& \sum_{n=0}^\infty  B_{2n+1} s^n \\ 
&=& 1 + \sum_{n=1}^\infty  \biggl[ {1 \over 2}\left[ \sum_{k=0}^{2n} B_{2n-k} B_k \right] + {1\over
  2} B_n \biggr] s^n\\ 
&=& 1 + {1 \over 2} \sum_{n=1}^\infty \sum_{k=0}^{n-1}
\biggl[B_{2n-2k-1}B_{2k+1}+ B_{2n-2k} B_{2k} \biggr]s^n + {1\over 2}
  \sum_{n=1}^\infty \left 
  ( B_0 B_{2n} + B_n \right) s^n\\
&=& 1 + {s\over 2} \sum_{n=0}^\infty \sum_{k=0}^\infty \left[ B_{2n+1} B_{2k+1}
  + B_{2n+2} B_{2k} \right] s^{n+k} + {1\over 2} \sum_{n=1}^\infty (B_{2n}
  + B_n) s^n\\
&=& 1+ {s \over 2} (\sum_{n=0}^\infty B_{2n+1} s^{n}) ( \sum_{k=0}^\infty
  B_{2k+1} s^k) + {1 \over 2} (\sum_{n=0}^\infty B_{2n+2} s^{n+1}) (\sum_{k=0}^\infty
  B_{2k} s^k)\\
&& + {1\over 2} (H(s)+G(s) - 2)\\
&=& {G(s)+ s I(s)^2 +H(s)^2 \over 2}
\end{eqnarray*}

\subsection{Summary of Generating Function Formulas}

\begin{eqnarray}
\nonumber H(s)&=& 1+ s I(s)H(s) \\
\nonumber I(s)&=& {1\over 2} \biggl[ G(s)+ s I(s)^2 +H(s)^2\biggr]\\
\nonumber G(s)&=& H(s^2) + s I(s^2)
\end{eqnarray}
The known \cite{Etherington,etherington1937non}  generating function is:
\begin{eqnarray}
\label{GFE} G(s) &=& 1 + {s \over 2}\biggl[G(s)^2 + G(s^2)\biggr] 
\end{eqnarray}

\section{Number of Non-Isomorphic Rooted Binary Trees with Number of Non-Isomorphic Siblings Given}
\label{NNIRBTWNNISG}

As discussed before, the functional equation \ref{GFE} is a 
generating function which counts the number of non-isomorphic trees with
labeled roots with $n$ nodes. 

Graph isomorphism is an equivalence relation on graphs and 
as such it partitions the class of all graphs into equivalence classes.  Therefore, an equivalence 
class is defined by the tree isomorphism with labeled roots:  
binary tree A with labeled root $a$ is isomorphic to binary
tree B with labeled root node $b$. This is clearly symmetric, reflexive and
transitive and is therefore an equivalence class. The functional equation 
\ref{GFE} above defines the number of binary trees
in that equivalence class based on the coefficient of $s$ to the number of nodes $n$.

It is interesting to develop a function which designates the cardinality of these
equivalence classes in terms of the number of ordered trees
therein and the
multiplicity of equivalence classes with the same cardinality. In order to
count these objects, observe the cardinality of these equivalence classes is a perfect 
power 
of two in all cases.  This follows from simply counting the number of siblings, which are the roots of non-isomorphic subtrees with labeled roots.  Each
instance of such a case implies two times the number 
of ordered trees which are isomorphic to it. This
parameter forms a variable which facilitates the development of a recurrence
for these objects.   Note that a complete tree has no
siblings, which are 
non-isomorphic in this sense.   At the other extreme, a tree of maximum depth
forms the maximum of this parameter.  Namely, siblings are all non-isomorphic due to a 
counting argument. 

With significant attention to all the sub cases, one may write the following
recurrence. $K_{n,\ell}$ is the number of equivalence classes of cardinality $2^\ell$.
Note that each member of the equivalence class consists of binary trees with
$n$ nodes and $\ell$ non-isomorphic sibling subtrees with labeled roots. 

\begin{eqnarray*}
K_{2n,\ell}&=& {1 \over 2} \sum_{k=0}^{2n-1} \sum_{v=0}^{\ell-1} K_{k,v}
K_{2n-1-k,\ell-1-v} \\
K_{2n+1,2\ell} &=& K_{n,\ell} + {1\over 2} \sum_{k=0}^{2n} \sum_{v=0}^{2\ell-1} K_{k,v}
K_{2n-k,2\ell-1-v} \\
K_{2n+1,2\ell+1} &=& {1 \over 2} \biggl(K_{n,\ell} \bigl(K_{n,\ell}-1\bigr)- K_{n,\ell}^2+
\sum_{k=0}^{2n} \sum_{v=0}^{2\ell} K_{k,v}
K_{2n-k,2\ell-v} \biggr) \\
&=&{1 \over 2} \biggl(- K_{n,\ell}+
\sum_{k=0}^{2n} \sum_{v=0}^{2\ell} K_{k,v}
K_{2n-k,2\ell-v} \biggr)
\end{eqnarray*}
 
The following base cases take precedence over the recurrence relations. 

\begin{eqnarray*}
K_{n,0}&=& \Biggl\{ {{1, \quad n=2^\ell -1, \quad \ell
\in \mathbb{Z}^{+,0}}\atop{0, \qquad \qquad \qquad \mathrm{otherwise}}}
\biggr.\label{bc0}\\ 
K_{s,t}&=& 0, \quad t \geq s > 0\label{bc2}\\
K_{s,s-1}&=& 1, \quad s > 1\label{bc3}
\end{eqnarray*}

\begin{figure}[h]
\begin{center}
\ifpdf
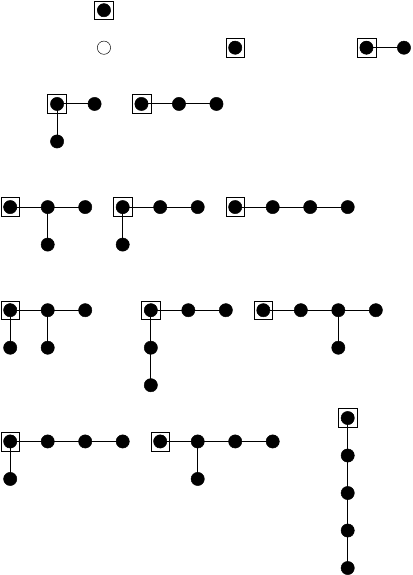
\else
\input{isonumnonisosiblings.pstex_t}
\fi
\end{center}
\caption{Rooted Binary Trees, Parameritized with Number of Nodes=$n$, and Number of Non-Isomorphic Siblings=$\ell$}
\label{rbtnell}
\end{figure}

Assume if $a$ and $b$ are outside of the region of definition
$K_{a,b}=0$. Table \ref{tbla} indicates the first few values for $K_{n,\ell}$.   
Figure \ref{rbtnell} illustrates the unique tree structure for the first few values. 

\begin{table}
\begin{center}
\begin{tabular}{|l|llllllllllll|}
\hline
\backslashbox{$n$}{$\ell$}&0&1&2&3&4&5&6&7&8&9&10&11\\
\hline
0&1&&&&&&&&&&&\\
1&1&0&&&&&&&&&&\\ 
2&0&1&0&&&&&&&&&\\
3&1&0&1&0&&&&&&&&\\
4&0&1&1&1&0&&&&&&&\\
5&0&1&2&2&1&0&&&&&&\\
6&0&0&3&3&4&1&0&&&&&\\
7&1&0&1&7&7&6&1&0&&&&\\
8&0&1&1&6&14&14&9&1&0&&&\\
9&0&1&3&4&21&28&28&12&1&0&&\\
  10& 0&0&3&8&17&54&58&50&16&1&0&\\
  11& 0&1&2&9&27&61&126&119&85&20&1&0\\
\hline
\end{tabular}
\hfill\\
\hfill\\
\end{center}
\caption{First few values of $K_{n,\ell}$}
\label{tbla}
\end{table}

Note that the following equations hold. This seems to indicate that a
relationship between the bivariate generating function for $K_{n,\ell}$ and the
aforementioned generating functions $F(z)$ and $G(s)$ may be developed.

\begin{eqnarray*}
C_n&=& \sum_{k=0}^{n} 2^k K_{n,k}\\
B_n&=& \sum_{k=0}^{n} K_{n,k}
\end{eqnarray*}

That is why the above recurrence is important to solving the functional
equation.  It forms a connection between the generating function $G(s)$ and
$R(z)$. 

One may modify the recursion by taking account of the base case equations \ref{bc0}, \ref{bc2} and
\ref{bc3} to obtain a more efficient summation of cases. 

\begin{eqnarray*}
K_{2n,\ell}&=& {1 \over 2} \sum_{k=0}^{2n-1}
\sum_{v=\max(0,\ell+k-2n)}^{\min(\ell-1,k)} K_{k,v} 
K_{2n-1-k,\ell-1-v} \\
K_{2n+1,2\ell} &=& K_{n,\ell} + {1\over 2} \sum_{k=0}^{2n}
\sum_{v=\max(0,k+2\ell-2n-1)}^{\min(2\ell-1,k)} K_{k,v} 
K_{2n-k,2\ell-1-v} \\
K_{2n+1,2\ell+1} &=& {1 \over 2} \biggl(- K_{n,\ell}+
\sum_{k=0}^{2n} \sum_{v=\max(0,k+2\ell-2n)}^{\min(2\ell,k)} K_{k,v}
K_{2n-k,2\ell-v} \biggr)
\end{eqnarray*}

%% moved to trove 2A
% recognized it's easier to work with the above than individual eqns.

\subsection{Bivariate Generating Function for $K_{n,\ell}$ Recurrence}

Define the bivariate generating function and solve the equation for
the functional equation for $ L (x,y)$.

\begin{eqnarray*}
 L (x,y)&=&\sum_{n=0}^\infty \sum_{\ell=0}^n
K_{n,\ell}\; x^n y^\ell 
\end{eqnarray*}

Define the following three utility functions to
express the system of functional equations. The functions mirror the
three cases in the recurrence relation for $K_{n,\ell}$. 

\begin{eqnarray*}
 P (x,y)&=&\sum_{n=0}^\infty \sum_{\ell=0}^{2n}
K_{2n,\ell}\; x^n y^\ell \\
 Q (x,y)&=&\sum_{n=0}^\infty \sum_{\ell=0}^{n}
K_{2n+1,2\ell}\; x^n y^\ell \\
 R (x,y)&=&\sum_{n=0}^\infty \sum_{\ell=0}^{n}
K_{2n+1,2\ell+1}\; x^n y^\ell 
\end{eqnarray*}

Also define the functions,\footnote{Henceforth, for compactness of notation
  all summations are assumed to step in increments of one.  If the upper limit
  of the summation is less than the lower limit of the summation, the entire term is
 taken to be zero.}

\begin{eqnarray*}
 P _0(x,y)&=&\sum_{n=0}^\infty \sum_{\ell=0}^{n}
K_{2n,2\ell}\; x^n y^\ell \\
 P _1(x,y)&=&\sum_{n=0}^\infty \sum_{\ell=0}^{n-1}
K_{2n,2\ell+1}\; x^n y^\ell \\
\end{eqnarray*}

Note that 

\[
 P (x,y)=  P _0(x,y^2) + y  P _1(x,y^2)
\]

and clearly

\[
 L (x,y)=  P (x^2,y)+
x\;  Q (x^2,y^2)+x y \; R (x^2,y^2) 
\]

\begin{eqnarray}
\label{Q}
 Q (x,y) =   L (x,y) + y  P _0(x,y)  P _1(x,y) + x
y  Q (x,y)  R (x,y)  
\end{eqnarray}

and thus

\[
 Q (x,y) =  {{ L (x,y) + y  P _0(x,y)  P _1(x,y)} \over
{1-x y  R (x,y)  }}
\]

A similar set of manipulations of the indices yields the functional equation.

\begin{eqnarray}
\label{R}
 R (x,y) = {1\over 2} \Bigl( - L (x,y)+ x
    \bigl(  Q (x,y)^2+ y 
     R (x,y)^2 \bigr)+ P _0(x,y)^2+ y  P _1(x,y)^2 \Bigr) 
\end{eqnarray}

Solving for $ R (x,y)$ in the quadratic one obtains.

\[
 R (x,y) = {{1 \pm \sqrt{1-x y \bigl(- L (x,y)+x
       Q (x,y)^2+ P _0(x,y)^2+y  P _1(x,y)^2  \bigr)
      }}\over{x y}}
\]

Similarly,

\begin{eqnarray}
\label{P0}
 P _0(x,y) = 1+ x y \bigl(  R (x,y)  P _0(x,y) +
 Q (x,y)  P _1(x,y) \bigr)
\end{eqnarray}

and 

\begin{eqnarray}
\label{P1}
 P _1(x,y) =  x  P _0(x,y)  Q (x,y) + x y
 P _1(x,y)  R (x,y)
\end{eqnarray}

\subsection{Summary of System of Functional Equations}

\begin{eqnarray*}
 L (x,y)&=&  P (x^2,y)+
x\;  Q (x^2,y^2)+x y \; R (x^2,y^2) \\
 P (x,y)&=&  P _0(x,y^2) + y  P _1(x,y^2)\\
 P (x,y)& = &{1\over{ 1- x y 
\left( Q (x,y^2) + y  R (x,y^2)\right)}}\\
 Q (x,y) &=&  {{ L (x,y) + y  P _0(x,y)  P _1(x,y)} \over
{1-x y  R (x,y)  }}\\
 R (x,y) &=& {{1 \pm \sqrt{1-x y \bigl(- L (x,y)+x
       Q (x,y)^2+ P _0(x,y)^2+y  P _1(x,y)^2  \bigr)
      }}\over{x y}}\\
 P _0(x,y) &=& {{1+  x y
 Q (x,y)  P _1(x,y)}\over{1 - x y  R (x,y)}} \\
 P _1(x,y) &=&  {{x  P _0(x,y)  Q (x,y)}\over{ 1- x y  R (x,y)}}
\end{eqnarray*}

Note that the latter two equations yield

\begin{eqnarray*}
 P _0(x,y) &=& {{1 - x y R (x,y)}\over{\bigl(1 - x y R (x,y)\bigr)^2-x^2  Q (x,y)^2}} \\
 P _1(x,y) &=& {{x Q (x,y)}\over{\bigl(1 - x y
 R (x,y)\bigr)^2-x^2  Q (x,y)^2}} \\
\end{eqnarray*}

Consider the functional equations \ref{P0}, \ref{P1}, \ref{Q} and \ref{R}. 

\begin{eqnarray*}
 Q (x,y) &=&   L (x,y) + y  P _0(x,y)  P _1(x,y) + x
y  Q (x,y)  R (x,y)  \\
 R (x,y) &=& {1\over 2} \Bigl( - L (x,y)+ x
    \bigl(  Q (x,y)^2+ y 
     R (x,y)^2 \bigr)+ P _0(x,y)^2+ y  P _1(x,y)^2 \Bigr) \\
 P _0(x,y) &=& 1+ x y \bigl(  R (x,y)  P _0(x,y) +
 Q (x,y)  P _1(x,y) \bigr)\\
 P _1(x,y) &=&  x  P _0(x,y)  Q (x,y) + x y
 P _1(x,y)  R (x,y)
\end{eqnarray*}

\begin{eqnarray*}
 L (x,y)&=&  P _0(x^2,y^2) + y  P _1(x^2,y^2)+
x\;  Q (x^2,y^2)+x y \; R (x^2,y^2) \\
 L (x,y)^2&=& P _0(x^2,y^2)^2 + y^2  P _1(x^2,y^2)^2\nonumber\\
&&{} +
x^2\;  Q (x^2,y^2)^2+x^2 y^2 \; R (x^2,y^2)^2\nonumber\\
&& {} + 2 y  P _0(x^2,y^2)  P _1(x^2,y^2) + 2 x P _0(x^2,y^2)  Q (x^2,y^2) \nonumber\\
&& {} + 2 x y  P _0(x^2,y^2)  R (x^2,y^2)\nonumber + 2 x y  P _1(x^2,y^2)  Q (x^2,y^2) \\
&& {} + 2 x y^2  P _1(x^2,y^2)  R (x^2,y^2) + 2 x^2 y \;  Q (x^2,y^2)  R (x^2,y^2) 
\end{eqnarray*}

The system can be reduced to a single functional equation, by matching the squared terms and substituting the cross terms of the squared generating function.  

\begin{eqnarray*}
L(x,y)^2&=& 2 R(x^2,y^2)+ L(x^2,y^2) + {2 \over y} \left(Q(x^2,y^2) - L(x^2,y^2) \right)\nonumber\\
&& {} + {2 \over x} P_1(x^2,y^2) + {2 \over x y} \left( P_0(x^2,y^2)-1\right)
\end{eqnarray*}

Which implies 

\begin{eqnarray*}
x y L(x,y)^2 &=&  2 \left( P_0(x^2,y^2) + y P_1(x^2,y^2) + x Q(x^2,y^2) + x y R(x^2,y^2)\right)\nonumber \\
&& {}  + x y L(x^2,y^2) - 2 x L(x^2,y^2) -2\\
&=& 2 L(x,y) +  x y L(x^2,y^2) - 2 x L(x^2,y^2) -2
\end{eqnarray*}

\begin{eqnarray*}
{x y L(x,y)^2\over 2} +1 &=&   L(x,y)+  x \left({y \over 2} -1\right) L(x^2,y^2)
\end{eqnarray*}

% VALIDATED again 8/18/2013

\[
L(x,y)=1 +   {\frac{x\,y\,{{L(x,y)}^2}}{2}} + x\,L({x^2},{y^2}) - 
   {\frac{x\,y\,L({x^2},{y^2})}{2}} 
\]

\begin{align}
\label{GEN1} L(x,y)=1 +   {\frac{x\,y\,{{L(x,y)}^2}}{2}} + x \left(1-{y\over 2}\right)\,L({x^2},{y^2}) 
\end{align}

One may see that in the case of $y=1$ it devolves to the known functional
equation\ref{GFE}.

\[
L(x,1)= 1 + {x \over 2} \left ( L(x,1)^2 + L(x^2,1) \right ) 
\]

In the case of $y=2$ it devolves to the known generating function equation \ref{GFE0}.

\[
L(x,2)= 1 + x\,{{L(x,2)}^2} 
\]

where

\[
L(x,2) = {\frac{1 \pm {\sqrt{1 - 4\,x}}}{2\,x}}
\]

\section{2-colored Rooted Non-isomorphic Binary Trees}
\label{TCRNIBT}

Consider a binary tree which has nodes colored zero or one.  Now define
isomorphism to be the standard bijection of edges and nodes \cite{bouritsas2022improving}; however, also
require that the node colors are matched as well. By way of definition, let the 
following program represent the 2-color isomorphism.  It
is similarly symmetric, reflexive and transitive and thus an equivalence
relation just like graph isomorphism. 

\begin{program}
\quad \mathrm{struct\;tree}\; \{
\quad \quad \mathrm{struct\;tree}\; *left;
\quad \quad \mathrm{struct\;tree}\; *right;
\quad \quad \mathrm{int}\;color; 
\quad \}\;\mathrm{tree};
\end{program}

\begin{program}
\quad{{\mathrm{int}}\;\mathrm{two\_color\_iso(tree}} *a,{\mathrm{tree}}*b) \{
\quad\quad{\mathrm{if}\;}((a==0)\;\land\;(b==0))\;\mathrm{return\;1};
\quad\quad{\mathrm{if}\;}((a==0)\; \lor \;(b==0))\;\mathrm{return\;0};
\quad\quad{\mathrm{if}\;}((a\rightarrow color)\neq(b\rightarrow color))\;\mathrm{return\;0};
\quad\quad{\mathrm{if}}\;\bigl(\;\;{\mathrm{two\_color\_iso}(a\rightarrow right,b\rightarrow right)}\;\bigr.
\quad\quad\quad\land\bigl.\;{\mathrm{two\_color\_iso}}(a\rightarrow left,b\rightarrow left)\bigr)\;\mathrm{return}\;1;
\quad\quad{\mathrm{if}}\;\bigl(\;\;{\mathrm{two\_color\_iso}(a\rightarrow right,b\rightarrow left)}\;\bigr.
\quad\quad\quad\land\bigl.\;{\mathrm{two\_color\_iso}}(a\rightarrow left,b\rightarrow right)\bigr)\;\mathrm{return}\;1;
\quad\quad{\mathrm{return\;0}};
\quad\}
\end{program}

\begin{figure}[h]
\begin{center}
\resizebox{!}{3.5in}{
\ifpdf
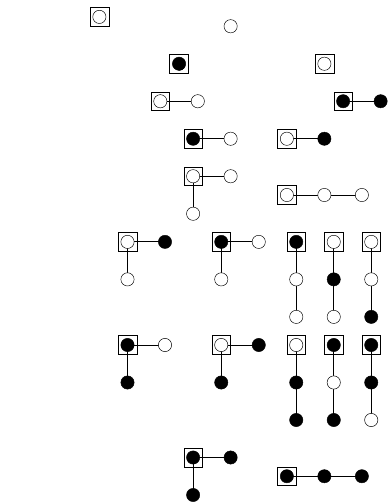
\else
\input{2colorParam.pstex_t}
\fi
}
\end{center}
\caption {2-Color Binary Trees, Parametrized by Number of Nodes=$n$ and Number of a Specific Color=$k$}
\label{twocolor}
\end{figure}

Given these 2-colored nodes in a binary tree structure, let $B_{n,k}$ be the number of such two colored rooted non-isomorphic binary trees with $n$ nodes and $k$ colored as black or $1$.  Table \ref{tbla2color} and Figure \ref{twocolor} show the first few examples.

\begin{align*}
B_{0,0}=B_{1,1}=B_{1,0}=&1\\
B_{0,k}=&0    & k>0\\
B_{1,k}=&0    &k>1\\
B_{n,0}=&B_n  &n>1\\
B_{n,k}=&0&     k>n
\end{align*}

Where $B_n$ is defined as above.

\[
B_0=B_1=1
\]

For $n\geq 1$, 

\begin{eqnarray*}
B_{2 n} &=&  {1 \over 2} \sum_{k=0}^{2n-1} B_{2n-k-1} B_k \\
B_{2 n+1} &=& {1 \over 2}\left(\left[ \sum_{k=0}^{2n} B_{2n-k} B_k \right] +
  B_n \right) 
\end{eqnarray*}

There are now three recursive cases to consider.  
For $n>0$ and $k>0$ this case can never have isomorphic subtrees due to a counting argument.

\begin{eqnarray*}
B_{2 n, k} &=&  {1 \over 2} \sum_{\ell=0}^{2n-1} \sum_{m=0}^{k-1} B_{\ell,m} B_{2n-1-\ell,k-1-m}+{1 \over 2} \sum_{\ell=0}^{2n-1} \sum_{m=0}^{k} B_{\ell,m} B_{2n-1-\ell,k-m}
\end{eqnarray*}

The first term corresponds to the case with the root node black; the second case corresponds to a white node as root. 

Now, similar to the above tree structure, counting an additional $B_{n,k}$ term shows up.

For $n\geq 0$ and $k>0$ we now have
\begin{eqnarray*}
B_{2n+1,2k} &=& {1 \over 2} \sum_{\ell=0}^{2n} \sum_{m=0}^{2k-1} B_{\ell,m} B_{2n-\ell,2k-1-m}\\
&&{}+{1 \over 2} \left (\left [\sum_{\ell=0}^{2n} \sum_{m=0}^{2k} B_{\ell,m} B_{2n-\ell,2k-m} \right] + B_{n,k} \right)
\end{eqnarray*}

The first term corresponds to the root node as black and the second term corresponds to the root node as white.  In the second term, the square bracketed term captures all of the possibilities; however, as before, with an excess of ${1\over 2} B_{n,k}$.  Adding an additional ${1\over 2} B_{n,k}$ fills in the possibilities where sibling subtrees are isomorphic \cite{yamazaki2024efficient}.  

For $n \geq 0$ and $k\geq 0$ we get
\begin{eqnarray*}
B_{2n+1,2k+1} &=& {1 \over 2} \left (\left [\sum_{\ell=0}^{2n} \sum_{m=0}^{2k} B_{\ell,m} B_{2n-\ell,2k-m} \right] + B_{n,k} \right)\\
&&{}+{1 \over 2} \sum_{\ell=0}^{2n} \sum_{m=0}^{2k+1} B_{\ell,m} B_{2n-\ell,2k+1-m}
\end{eqnarray*}

\begin{table}[h]
\begin{center}
\resizebox{!}{1in}{
\begin{tabular}{|l|lllllllllll|}
\hline
\backslashbox{$n$}{$k$}&0&1&2&3&4&5&6&7&8&9&10\\
\hline
0&1&&&&&&&&&&\\
1&1&1&&&&&&&&&\\
2&1&2&1&&&&&&&&\\
3&2&5&5&2&&&&&&&\\
4&3&11&16&11&3&&&&&&\\
5&6&26&50&50&26&6&&&&&\\
6&11&60&143&188&143&60&11&&&&\\
7&23&142&404&656&656&404&142&23&&&\\
8&46&334&1105&2143&2652&2143&1105&334&46&&\\
9&98&794&2995&6737&9934&9934&6737&2995&794&98&\\
10&207&1888&7999&20504&35080&41788&35080&20504&7999&1888&207\\
\hline
\end{tabular}}
\hfill\\
\hfill\\
\end{center}
\caption {Number of 2-Color Binary Trees, Parametrized by Number of Nodes=$n$ and Number of a Specific Color=$k$}
\label{tbla2color}
\end{table}

A similar explanation to the previous case is applicable.   
In the following, let summations always increment by one or not be summed.  Removing the need for the condition $B_{n,k}=0$ for $k>n$ we now obtain:

For $n>0$ and $k>0$ 
\begin{eqnarray*}
B_{2 n, k} &=&  {1 \over 2} \sum_{\ell=0}^{2n-1} \sum_{m=\max(0,\ell+k-2n)}^{\min(\ell,k-1)} B_{\ell,m} B_{2n-1-\ell,k-1-m}\\
&&{} +{1 \over 2} \sum_{\ell=0}^{2n-1} \sum_{m=\max(0,\ell+k+1-2n)}^{\min(\ell,k)} B_{\ell,m} B_{2n-1-\ell,k-m}\\
&=&  {1 \over 2} \sum_{r=0}^1 \sum_{\ell=0}^{2n-1} \sum_{m=\max(0,\ell+k+1-r-2n)}^{\min(\ell,k-r)} B_{\ell,m} B_{2n-1-\ell,k-r-m}\\
\end{eqnarray*}

Now for $n\geq 0$ and $k>0$ 
\begin{eqnarray*}
B_{2n+1,2k} &=& {1 \over 2} \sum_{\ell=0}^{2n} \sum_{m=\max(0,\ell+2k-1-2n)}^{\min(\ell,2k-1)} B_{\ell,m} B_{2n-\ell,2k-1-m}\\
&&{}+{1 \over 2} \left (\left [\sum_{\ell=0}^{2n} \sum_{m=\max(0,\ell+2k-2n)}^{\min(\ell,2k)} B_{\ell,m} B_{2n-\ell,2k-m} \right] + B_{n,k} \right)\\
&=& {1 \over 2} \sum_{r=0}^1 \sum_{\ell=0}^{2n} \sum_{m=\max(0,\ell+2k-r-2n)}^{\min(\ell,2k-r)} B_{\ell,m} B_{2n-\ell,2k-r-m}+{1 \over 2}  B_{n,k} 
\end{eqnarray*}

With $n \geq 0$ and $k\geq 0$ we get
\begin{eqnarray*}
B_{2n+1,2k+1} &=& {1 \over 2} \left (\left [\sum_{\ell=0}^{2n} \sum_{m=\max(0,\ell+2k-2n)}^{\min(\ell,2k)} B_{\ell,m} B_{2n-\ell,2k-m} \right] + B_{n,k} \right)\\
&&{} +{1 \over 2} \sum_{\ell=0}^{2n} \sum_{m=\max(0,\ell+2k+1-2n)}^{\min(\ell,2k+1)} B_{\ell,m} B_{2n-\ell,2k+1-m}\\
&=&{1 \over 2} \sum_{r=0}^1 \sum_{\ell=0}^{2n} \sum_{m=\max(0,\ell+2k+r-2n)}^{\min(\ell,2k+r)} B_{\ell,m} B_{2n-\ell,2k+r-m}  + {1 \over 2} B_{n,k} 
\end{eqnarray*}

Define the generating function $M(x,y)$.

\begin{eqnarray*}
M(x,y)&=&\sum_{n=0}^\infty \sum_{k=0}^n  B_{n,k} x^n y^k\\
&=&1+\sum_{n=1}^\infty \sum_{k=1}^n B_{n,k} x^n y^k + \sum_{n=1}^\infty B_{n,0} x^n \\
\end{eqnarray*}

Define the following useful generating functions.

\begin{eqnarray*}
P(x,y)&=& \sum_{n=0}^\infty \sum_{k=0}^{2n}  B_{2n,k} x^n y^k\\
Q(x,y)&=& \sum_{n=0}^\infty \sum_{k=0}^n B_{2n+1,2k} x^n y^k\\
R(x,y)&=& \sum_{n=0}^\infty \sum_{k=0}^n B_{2n+1,2k+1} x^n y^k\\
\end{eqnarray*}

Defining the auxiliary functions

\begin{eqnarray*}
P_0(x,y)&=& \sum_{n=0}^\infty \sum_{k=0}^n  B_{2n,2k} x^n y^k\\
P_1(x,y)&=& \sum_{n=1}^\infty \sum_{k=0}^{n-1} B_{2n, 2k+1} x^n y^k
\end{eqnarray*}

\[
P(x,y)= P_0(x,y^2) + y P_1(x,y^2)
\]

\[
M(x,y)= P(x^2,y)+x Q(x^2,y^2)+x y R(x^2,y^2) 
\]

% VALID Expression08x11x2013  8/11/2013
\begin{eqnarray*}
R(x,y)-{1 \over 2} M(x,y) &=& {1 \over 2}  \Biggl(  P_0(x,y)^2 + 2   P_0(x,y)P_1(x,y)\\
&& {}+ y P_1(x,y)^2 + x Q(x,y)^2 +x y R(x,y)^2 \\
&& \Biggl. {}+ 2   Q(x,y) R(x,y) \Biggr) \\
\end{eqnarray*}

So all four generating function equations are:

\begin{eqnarray}
\label{EQP0} P_0(x,y)&=&
%
%r=0, q=0
%  
x P_0(x,y) Q(x,y) + x y P_1(x,y)Q(x,y)-x P_0(x,0) Q(x,0) \nonumber \\
&& {} + x y P_0(x,y) R(x,y)+ x y P_1(x,y) R(x,y) + P_0(x,0) \\
\label{EQP1} P_1(x,y)&=& x P_0(x,y) Q(x,y)+ x P_1(x,y) Q(x,y)+ x P_0(x,y) R(x,y) \nonumber\\
&& {} + x y P_1(x,y) R(x,y)\\
\label{EQQ}
Q(x,y)&=& y P_0(x,y) P_1(x,y) + x y Q(x,y) R(x,y)+Q(x,0) \nonumber \\
&& {} + {1 \over 2} \Bigl( M(x,y) -M(x,0) + P_0(x,y) ^2  - P_0(x,0)^2 \Bigr.\nonumber\\
&& \Bigl.{} + y P_1(x,y)^2 + x Q(x,y)^2 - x Q(x,0)^2+ xy R(x,y)^2 \Bigr)\\
\label{EQR}
R(x,y) &=&    P_0(x,y)P_1(x,y) +  x Q(x,y) R(x,y)  + {1 \over 2}  \Bigl(  P_0(x,y)^2  \Bigr.\nonumber\\
&& \Bigl. {} +  y P_1(x,y)^2 + x Q(x,y)^2 +x y R(x,y)^2 +M(x,y)\Bigr) 
\end{eqnarray}

Expressing for $Q(x,y)$ and $R(x,y)$ as a quadratic from equations \ref{EQQ} and \ref{EQR} respectively.
% solved in mathematica then expressed a a quadratic.

\begin{flalign}
&\label{EQQQUAD} 0 ={x \over 2} [Q(x,y)]^2 - (1 - x y R(x,y)) [Q(x,y)] \nonumber &\\
& \quad {} + (M(x,y)-M(x,0)+P_0(x,y)^2-P_0(x,0)^2+2y P_0(x,y) P_1(x,y)\nonumber &\\
& \quad {} + y P_1(x,y)^2 +2 Q(x,0)-x Q(x,0)^2+ x y R(x,y)^2) &\\
&\label{EQRQUAD} 0= {x y\over 2} [R(x,y)]^2 - (1 - x Q(x,y)) [ R(x,y)]  \nonumber &\\
& \quad {} + (M(x,y) +P_0(x,y)^2 + 2 P_0(x,y) P_1(x,y) + y P_1(x,y)^2 + x Q(x,y)^2)&
\end{flalign}

Note also that following from the definitions:

\begin{eqnarray}
\label{MDEF} M(x,y)= P_0(x^2,y^2) + y P_1(x^2,y^2) +x Q(x^2,y^2)+x y R(x^2,y^2) 
\end{eqnarray}

For notation purposes only assume that $0^0=1$ (this can also be handled with limits) when expressed as $Q(x,0)$ for example in the following.  Solving equations \ref{EQP0} and \ref{EQP1} by first solving for $P_0(x,y)$ and $P_1(x,y)$ respectively then substituting into each equation.  

\begin{eqnarray}
\label{EQP0A} P_0(x,y) &=& {{P_0(x,0) (1-x Q(x,0))(1-x Q(x,y)-x y R(x,y))}\over
 {(1-x Q(x,y)-x y R(x,y))^2-x^2 y (Q(x,y)+R(x,y))^2}}\\
\label{EQP1A} P_1(x,y) &=& {{x P_0(x,0) (1-xQ(x,0))(Q(x,y)+R(x,y))}\over{(1- x Q(x,y)- x y R(x,y))^2- x^2 y (Q(x,y)+R(x,y))^2}}
\end{eqnarray}

Solving the quadratic in equations \ref{EQQQUAD} and \ref{EQRQUAD}.

\begin{eqnarray}
\label{QA} Q(x,y)&=& {1\over {2 x}} \left[ (1- x y R(x,y)) \pm \left[  (1- x y R(x,y))^2 -4 x (M(x,y)-M(x,0)\right.\right.\nonumber\\
&& \left. \left. {} +P_0(x,y)^2-P_0(x,0)^2+2y P_0(x,y) P_1(x,y)\right.\right.\nonumber\\
&& \left.\left.  {} + y P_1(x,y)^2 +2 Q(x,0)-x Q(x,0)^2+ x y R(x,y)^2) \right] ^{1\over 2}\right]\\ 
\label{RA} R(x,y)&=& {1 \over {2 x y}} \left[ (1 - x Q(x,y)) \pm \left[(1 - x Q(x,y))^2- 4 x y   (M(x,y) +P_0(x,y)^2 \right.\right.\nonumber\\ 
&&\left.\left. {}+ 2 P_0(x,y) P_1(x,y) + y P_1(x,y)^2 + x Q(x,y)^2)\right]^{1 \over 2} \right]
\end{eqnarray}

Equations \ref{QA} and \ref{RA} imply that respectively:

\begin{flalign}
&M(x,y) = M(x,0)-P_0(x,y)^2+P_0(x,0)^2 - 2y P_0(x,y) P_1(x,y) - y  P_1(x,y)^2 & \nonumber \\
&{}+2Q(x,y) - x Q(x,y)^2 -2 Q(x,0)+ xQ(x,0)^2-2 x y Q(x,y) R(x,y)- x y R(x,y)^2&\\
&M(x,y)= -P_0(x,y)^2 - 2 P_0(x,y) P_1(x,y)- y P_1(x,y)^2-x Q(x,y)^2 &\nonumber\\
&{}+2 R(x,y)- 2 x Q(x,y) R(x,y) - x y  R(x,y)^2 &
\end{flalign}

Summarizing

\begin{flalign}
&M(x,y)= P_0(x^2,y^2) + y P_1(x^2,y^2) +x Q(x^2,y^2)+x y R(x^2,y^2) & \\
&\label{M2} M(x,y) = M(x,0)-P_0(x,y)^2+P_0(x,0)^2 - 2y P_0(x,y) P_1(x,y) &\nonumber \\
&\quad\quad{} - y  P_1(x,y)^2+2Q(x,y) - x Q(x,y)^2 -2 Q(x,0) &\\
&\quad\quad{}+  xQ(x,0)^2-2 x y Q(x,y) R(x,y)- x y R(x,y)^2 &\\
&M(x,y)= -P_0(x,y)^2 - 2 P_0(x,y) P_1(x,y)- y P_1(x,y)^2-x Q(x,y)^2&\nonumber\\
&\quad\quad {}+2 R(x,y)- 2 x Q(x,y) R(x,y) - x y  R(x,y)^2&\\
&P_0(x,y) = {{P_0(x,0) (1-x Q(x,0))(1-x Q(x,y)-x y R(x,y))}\over
 {(1-x Q(x,y)-x y R(x,y))^2-x^2 y (Q(x,y)+R(x,y))^2}}&\\
&P_1(x,y) = {{x P_0(x,0) (1-xQ(x,0))(Q(x,y)+R(x,y))}\over{(1- x Q(x,y)- x y R(x,y))^2- x^2 y (Q(x,y)+R(x,y))^2}}&
\end{flalign}

Note that Section \ref{SIMPLEGEN} equations are

\begin{eqnarray}
\label{H2} H(s)&=& 1+ s I(s)H(s) \\
\label{I2}I(s)&=& {1\over 2} \biggl[ G(s)+ s I(s)^2 +H(s)^2\biggr]\\
G(s)&=& H(s^2) + s I(s^2)\\
G(s) &=& 1 + {s \over 2}\biggl[G(s)^2 + G(s^2)\biggr] 
\end{eqnarray}

These equate to $H(x)= P_0(x,0)$ and $I(x)= Q(x,0)$.   $G(x)= M(x,0)$ therefore using equation \ref{I2} in equation \ref{M2} yields

\begin{eqnarray}
\label{M3} M(x,y) &=& -P_0(x,y)^2- 2y P_0(x,y) P_1(x,y) - y  P_1(x,y)^2+2Q(x,y)\nonumber \\
&&{}- x Q(x,y)^2 -2 x y Q(x,y) R(x,y)- x y R(x,y)^2
\end{eqnarray}

Equation \ref{H2} implies that $P_0(x,0) -x P_0(x,0) Q(x,0) = P_0(x,0) + (1- P_0(x,0))= 1$.  This and equation \ref{M3} simplify the system of equations to

\begin{eqnarray}
M(x,y)&=& P_0(x^2,y^2) + y P_1(x^2,y^2) +x Q(x^2,y^2)+x y R(x^2,y^2) \\
\label{MFirst} M(x,y) &=& 2 Q(x,y)-P_0(x,y)^2- 2y P_0(x,y) P_1(x,y) - y  P_1(x,y)^2\nonumber \\
&&{}- x Q(x,y)^2 -2 x y Q(x,y) R(x,y)- x y R(x,y)^2\\
\label{MSecond} M(x,y)&=& 2 R(x,y)-P_0(x,y)^2 - 2 P_0(x,y) P_1(x,y)- y P_1(x,y)^2\nonumber\\
&&{} -x Q(x,y)^2- 2 x Q(x,y) R(x,y) - x y  R(x,y)^2\\
P_0(x,y) &=& {{(1-x Q(x,y)-x y R(x,y))}\over
 {(1-x Q(x,y)-x y R(x,y))^2-x^2 y (Q(x,y)+R(x,y))^2}}\\
P_1(x,y) &=& {{x (Q(x,y)+R(x,y))}\over{(1- x Q(x,y)- x y R(x,y))^2- x^2 y (Q(x,y)+R(x,y))^2}}
\end{eqnarray}

Equating equations \ref{MFirst} and \ref{MSecond} yields

\begin{eqnarray}
\label{EQUATION40} Q(x,y) -R(x,y) &=&  (y-1)  (P_0(x,y) P_1(x,y) + x Q(x,y) R(x,y)) 
\end{eqnarray}

Substituting back in to equations \ref{MFirst} and \ref{MSecond} eliminating the cross terms and obtaining

\begin{eqnarray}
\label{MFirstOnlySquares} M(x,y) &=& 2 Q(x,y)-P_0(x,y)^2  + 2 y  \left({ {Q(x,y)-R(x,y)} \over {1-y} }\right) - y  P_1(x,y)^2\nonumber \\
&&{}- x Q(x,y)^2 - x y R(x,y)^2\\
\label{MSecondOnlySquares} M(x,y)&=& 2 R(x,y)-P_0(x,y)^2 + 2  \left({ {Q(x,y)-R(x,y)} \over {1-y} }\right)- y P_1(x,y)^2\nonumber\\
&&{} -x Q(x,y)^2 - x y  R(x,y)^2
\end{eqnarray}

Simplifying yields one equation for $M(x,y)$

% valid 8/16/2013
\begin{flalign}
&\label{MEQUATE} M(x,y)= 2\left( {Q(x,y) - y R(x,y) \over {1-y}}\right)-P_0(x,y)^2   - y  P_1(x,y)^2\nonumber &\\
&\quad\quad {}- x Q(x,y)^2 - x y R(x,y)^2 &\\
& (1-y)M(x,y)= 2\left( Q(x,y) - y R(x,y) \right)\nonumber &\\
& \quad\quad {}- (1-y)\left(P_0(x,y)^2   + y  P_1(x,y)^2+ x Q(x,y)^2 + x y R(x,y)^2\right)&
\end{flalign}

Going back to initial work on $P_0(x,y)$ and $P_1(x,y)$

\begin{flalign}
&\label{EQP0B} P_0(x,y)=&\nonumber \\
%
%r=0, q=0
%  
&\quad \quad x P_0(x,y) Q(x,y) + x y P_1(x,y)Q(x,y) &\nonumber \\
& \quad\quad {} + x y P_0(x,y) R(x,y)+ x y P_1(x,y) R(x,y) + 1 &\\
&\label{EQP1B} P_1(x,y) =& \nonumber \\
& \quad\quad x P_0(x,y) Q(x,y)+ x P_1(x,y) Q(x,y)\nonumber &\\
& \quad\quad {} + x P_0(x,y) R(x,y) + x y P_1(x,y) R(x,y)&
\end{flalign}

Solving for $P_0(x,y)$ and $P_1(x,y)$ without moving to just dependence on $Q(x,y)$ and $R(x,y)$ (that previous calculation gave a hint for subsitutions)

\begin{eqnarray}
P_0(x,y) &=& {1 + x y P_1(x,y) Q(x,y) + x y P_1(x,y) R(x,y) \over {1- x Q(x,y)- x y R(x,y)}}\\
P_1(x,y) &=&{ x P_0(x,y) Q(x,y) + x P_0(x,y) R(x,y) \over {1- x Q(x,y)- x y R(x,y)}}
\end{eqnarray}

Summarizing again 

\begin{flalign}
&\label{EQM11}M(x,y)=  P_0(x^2,y^2) + y P_1(x^2,y^2) +x Q(x^2,y^2)+x y R(x^2,y^2) & \\
&\label{EQM12}(1-y)M(x,y)= 2\left( Q(x,y) - y R(x,y) \right)&\nonumber\\
&\quad\quad{}- (1-y)\left(P_0(x,y)^2   + y  P_1(x,y)^2+ x Q(x,y)^2 + x y R(x,y)^2\right)&\\
&P_0(x,y) = {1 + x y P_1(x,y) Q(x,y) + x y P_1(x,y) R(x,y) \over {1- x Q(x,y)- x y R(x,y)}}&\\
&P_1(x,y) ={ x P_0(x,y) Q(x,y) + x P_0(x,y) R(x,y) \over {1- x Q(x,y)- x y R(x,y)}} &
\end{flalign}

Squaring in equation \ref{EQM11}, $M(x,y)$ yields the following and then substitution can occur.

\begin{flalign}
&\label{EQM13} M(x,y)^2= P_0(x^2,y^2)^2 + y^2 P_1(x^2,y^2)^2 + x^2 Q(x^2,y^2)^2 + x^2 y^2 R(x^2,y^2)^2&\nonumber\\
& \quad {} + 2 y P_0(x^2,y^2) P_1(x^2,y^2) + 2 x P_0(x^2,y^2) Q(x^2,y^2) + 2 x y P_1(x^2,y^2) Q(x^2,y^2)&\nonumber\\
& \quad {} + 2 x y P_0(x^2,y^2) R(x^2,y^2) + 2 x y^2 P_1(x^2,y^2) R(x^2,y^2) + 2 x^2y Q(x^2,y^2) R(x^2,y^2) &
\end{flalign}

Restating equations \ref{MFirst} and \ref{MSecond}

\begin{eqnarray*}
M(x,y) &=& 2 Q(x,y)-P_0(x,y)^2- 2y P_0(x,y) P_1(x,y) - y  P_1(x,y)^2\nonumber \\
&&{}- x Q(x,y)^2 -2 x y Q(x,y) R(x,y)- x y R(x,y)^2\\
M(x,y)&=& 2 R(x,y)-P_0(x,y)^2 - 2 P_0(x,y) P_1(x,y)- y P_1(x,y)^2\nonumber\\
&&{} -x Q(x,y)^2- 2 x Q(x,y) R(x,y) - x y  R(x,y)^2
\end{eqnarray*}

Taking the non squared parts from \ref{MFirst} and \ref{MSecond} and equating them to \ref{MEQUATE}

\begin{flalign}
&\label{QNOW} 2\left( {Q(x,y) - y R(x,y) \over {1-y}}\right)=\nonumber&\\
&\quad\quad\quad 2 Q(x,y)- 2y P_0(x,y) P_1(x,y) -2 x y Q(x,y) R(x,y)&\\
&\quad{}=2 R(x,y) - 2 P_0(x,y) P_1(x,y) - 2 x Q(x,y) R(x,y)&
\end{flalign}

Solving for the cross terms only.

\begin{eqnarray*}
2\left( {Q(x,y)  -  R(x,y)\over {1-y}}\right)&=&  - 2 P_0(x,y) P_1(x,y) - 2 x Q(x,y) R(x,y)
\end{eqnarray*}

Note that the cross terms in equation \ref{EQM13} are equal to the following

\begin{flalign*}
&2 y \left( {Q(x^2,y^2)  -  R(x^2,y^2)\over {1-y^2}}\right)=& \\
& \quad\quad {} - 2 y  P_0(x^2,y^2) P_1(x^2,y^2) - 2 x^2 y Q(x^2,y^2) R(x^2,y^2)&\\
\end{flalign*}

Using \ref{EQP0B} and \ref{EQP1B} and subtracting equations and subtracting after multiplication by $y$ as well.  This gives cross terms in the squared equation \ref{EQM13}.

\begin{eqnarray*}
P_1(x,y)-P_0(x,y)&=& (1-y) (x P_1(x,y) Q(x,y)   + x P_0(x,y) R(x,y)) -1\\
P_0(x,y)-y P_1(x,y) &=& (1-y) ( x P_0(x,y) Q(x,y) + x y P_1(x,y) R(x,y)) +1 
\end{eqnarray*}

\begin{eqnarray*}
{P_1(x,y)-P_0(x,y)+1 \over x( 1-y) }&=&  P_1(x,y) Q(x,y)   +  P_0(x,y) R(x,y)\\
{P_0(x,y)-y P_1(x,y)-1 \over x( 1-y)} &=&   P_0(x,y) Q(x,y) +  y P_1(x,y) R(x,y)
\end{eqnarray*}

Note that the cross terms in equation \ref{EQM13} are equal to the following:

\begin{flalign*}
& 2 x y \left({P_1(x^2,y^2)-P_0(x^2,y^2)+1 \over x^2( 1-y^2) }\right)=&\\
& \quad \quad 2 x y P_1(x^2,y^2) Q(x^2,y^2)   +  2 x y P_0(x^2,y^2) R(x^2,y^2) &\\
&2 x \left({P_0(x^2,y^2)-y^2 P_1(x^2,y^2)-1 \over x^2( 1-y^2)}\right) = & \\
& \quad\quad 2 x  P_0(x^2,y^2) Q(x^2,y^2) + 2 x y^2 P_1(x^2,y^2) R(x^2,y^2) &
\end{flalign*}

This is not factorization, but a removal of the cross terms and using \ref{MEQUATE} in equation \ref{EQM13} to eliminate the square terms yields:

%validated 8/18/13
\begin{eqnarray*}
M(x,y)^2 &=&  2\left( {Q(x^2,y^2) - y^2 R(x^2,y^2) \over {1-y^2}}\right) - M(x^2,y^2)\nonumber\\
&& {} - 2 y \left( {Q(x^2,y^2)  -  R(x^2,y^2)\over {1-y^2}}\right)\nonumber\\
&& {} + 2 x y \left({P_1(x^2,y^2)-P_0(x^2,y^2)+1 \over x^2( 1-y^2) }\right)\nonumber\\
&& {} + 2 x \left({P_0(x^2,y^2)-y^2 P_1(x^2,y^2)-1 \over x^2( 1-y^2)}\right)
\end{eqnarray*}

\begin{flalign*}
& x^2 (1-y^2)  \left (   M(x,y)^2 + M(x^2,y^2)  \right)  = &\\
 &\quad \quad2 x \left( -1 + P_0(x^2,y^2) +y P_1(x^2,y^2) +x Q(x^2,y^2) + x y R(x^2,y^2) \right)  & \\
 &\quad \quad {} - 2 x y \left ( -1 + P_0(x^2,y^2) + y P_1( x^2, y^2) + x Q(x^2,y^2) + x y R(x^2,y^2) \right )&
\end{flalign*}

\begin{eqnarray*}
x^2 (1-y^2)  \left (   M(x,y)^2 + M(x^2,y^2)  \right) &=& 2 x  (1-y) \left( M(x,y)-1 \right) 
\end{eqnarray*}

\begin{eqnarray*}
x (1+y)  \left (   M(x,y)^2 + M(x^2,y^2)  \right) &=& 2  \left( M(x,y)-1 \right) 
\end{eqnarray*}

The final generating function equation is 

% valid 08/18/13
\begin{eqnarray}
 \label{GEN2} M(x,y) &=& x (1+y)  \left (  { M(x,y)^2 + M(x^2,y^2)\over 2}  \right) +1
\end{eqnarray}

Which reduces to the known generating function \ref{GFE} when $y=0$ and $0^0$ is formally interpereted as $1$ (this can also be handled with limits)

\begin{eqnarray*}
 M(x,0) &=& x \left (  { M(x,0)^2 + M(x^2,0)\over 2}  \right) +1
\end{eqnarray*}

\section{2-Color Rooted Binary Tree Isomorphism, Parameritized with Number of Specific Color, Non-Isomorphic Siblings and Nodes}
\label{TCRBTIPNSCNISN}

One may now define a recurrence which can calculate the multiplicity of
equivalence classes of cardinality $2^\ell$: $K_{n,\ell,c}$ -- where $n$ is the
number of 
nodes and $c$ the number of nodes colored black or $1$; note that $\ell$ is the number of
non-isomorphic (under the isomorphism defined above or identically ``flip equivalence'') sibling subtrees in the
tree.  $p,q,r \in \{0,1\}$.  Figures \ref{fig:rbX} and \ref{fig:rbY} show the first few trees and Table \ref{tbla2} shows the first few enumerated by the parameters. 

% resume this problem 8/18/2013

\begin{flalign*}
K_{2n+p,2\ell+q,2c+r} =
p\left(-\frac{1}{2} \right)^q K_{n,\ell,c} &\\
{} + {1 \over 2}\sum_{\delta=0}^1\;\;
\sum_{k=0}^{2n+p-1} \;\; \sum_{v=0}^{2\ell+q-1} \;\;\sum_{
\begin{subarray}{c}
m=0\\
2c+r \geq \delta
\end{subarray}
}^{2c+r-\delta}
&K_{k,v,m} K_{2n+p-1-k,2\ell+q-1-v,2c+r-\delta-m} &
\end{flalign*}
 
The following base cases take precedence over the recurrence relations.  

\begin{eqnarray*}
K_{0,0,0}&=&1\\
K_{n,\ell,c}&=& 0, \quad \ell \geq n > 0 \;\lor\; \ell < 0 \; \lor\; c<0  \\
K_{n,\ell,c}&=& \biggl( {n \atop c} \biggr), \quad \ell =
n-1 \;\land\; 0 \leq c \leq n  \\ 
K_{n,\ell,c}&=& 1, \quad r \in \mathbb{N}_1 \;\land\; n=2^{r}-1 \;\land\; \ell=0
\;\land\; 0 \leq c \leq n \\
K_{n,\ell,c}&=& 0, \quad r \in \mathbb{N}_1 \;\land\; n\neq 2^{r}-1 \;\land\; \ell=0 \\
\end{eqnarray*}

Where $\mathbb{N}_1=\{1,2,3,\ldots \}$.

\begin{table}[h]
\begin{center}
\resizebox{!}{3.5in}{
\begin{tabular}{l|lllllllll|l}
\backslashbox{$n$}{$\ell$}&0&1&2&3&4&5&6&7&8&$c$\\
\hline
0&1&&&&&&&&&0\\
\hline  
1&1&0&&&&&&&&0\\
&1&0&&&&&&&&1\\
\hline
2&0&1&0&&&&&&&0\\
&0&2&0&&&&&&&1\\
&0&1&0&&&&&&&2\\
\hline
3&1&0&1&0&&&&&&0\\
&1&1&3&0&&&&&&1\\
&1&1&3&0&&&&&&2\\
&1&0&1&0&&&&&&3\\
\hline  
4&0&1&1&1&0&&&&&0\\
&0&2&5&4&0&&&&&1\\
&0&2&8&6&0&&&&&2\\
&0&2&5&4&0&&&&&3\\
&0&1&1&1&0&&&&&4\\
\hline
5&0&1&2&2&1&0&&&&0\\
&0&3&5&13&5&0&&&&1\\
&0&4&9&27&10&0&&&&2\\
&0&4&9&27&10&0&&&&3\\
&0&3&5&13&5&0&&&&4\\
&0&1&2&2&1&0&&&&5\\
\hline
6&0&0&3&3&4&1&0&&&0\\
&0&0&12&15&27&6&0&&&1\\
&0&0&21&37&70&15&0&&&2\\
&0&0&24&50&94&20&0&&&3\\
&0&0&21&37&70&15&0&&&4\\
&0&0&12&15&27&6&0&&&5\\
&0&0&3&3&4&1&0&&&6\\
\hline
7&1&0&1&7&7&6&1&0&&0\\
&1&1&6&34&45&48&7&0&&1\\
&1&2&15&76&141&148&21&0&&2\\
&1&3&20&108&239&250&35&0&&3\\
&1&3&20&108&239&250&35&0&&4\\
&1&2&15&76&141&148&21&0&&5\\
&1&1&6&34&45&48&7&0&&6\\
&1&0&1&7&7&6&1&0&&7\\
\hline
8&0&1&1&6&14&14&9&1&0&0\\
&0&2&5&39&86&116&78&8&0&1\\
&0&2&11&109&249&426&280&28&0&2\\
&0&2&17&179&447&876&566&56&0&3\\
&0&2&20&206&540&1104&710&70&0&4\\
&0&2&17&179&447&876&566&56&0&5\\
&0&2&11&109&249&426&280&28&0&6\\
&0&2&5&39&86&116&78&8&0&7\\
&0&1&1&6&14&14&9&1&0&8\\
\hline
\end{tabular}}
\hfill\\
\hfill\\
\end{center}
\caption{Number of Rooted Binary Trees, Parameterized with Number of Nodes=$n$, Number of Non-Isomorphic Siblings=$\ell$ and Number of Specific Color=$c$}
\label{tbla2}
\end{table}

\begin{figure}[h]
\begin{center}
\ifpdf
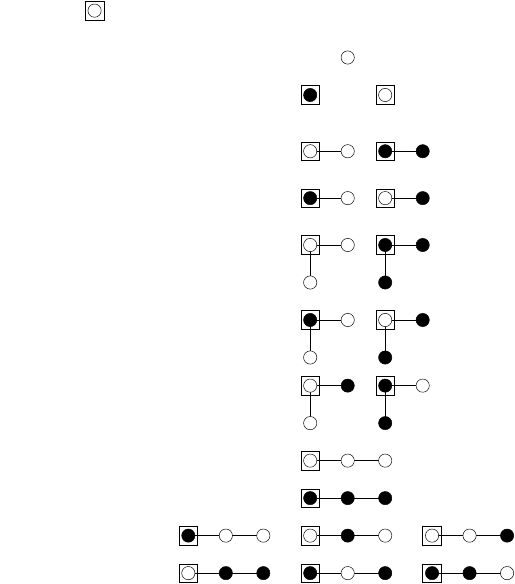
\else
\input{3ParamBTree.pstex_t}
\fi
\end{center}
\caption{Rooted Binary Trees, Parameritized with Number of Nodes=$n$, Number of Non-Isomorphic Siblings=$\ell$ and Number of Specific Color=$c$}
\label{fig:rbX}
\end{figure}

\begin{figure}[h]
\begin{center}
\ifpdf
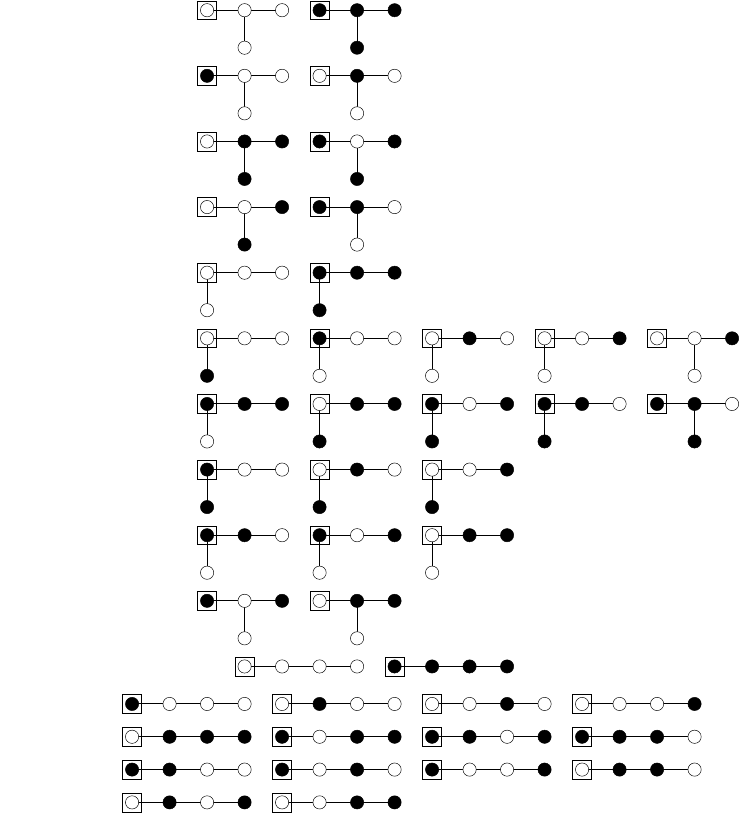
\else
\input{3ParamBTree2.pstex_t}
\fi
\end{center}
\caption{Continued:  Rooted Binary Trees, Parameritized with Number of Nodes=$n$, Number of Non-Isomorphic Siblings=$\ell$ and Number of Specific Color=$c$}
\label{fig:rbY}
\end{figure}

Define the following generating functions based on the $K_{n,\ell,c}$
recurrence. 

\begin{eqnarray}
\label{defnS}  S (x,y,z)&=&\sum_{n=0}^\infty \;\;\sum_{\ell=0}^n \sum_{c=0}^n K_{n,\ell,c} \;
x^n y^\ell z^c\\
\end{eqnarray}

\begin{align*}
 S _{p,q,r}(x,y,z)&=&\sum_{n=(1-p) (q+r- q r)}^\infty \sum_{\ell=0}^{n-(1-p)q}
\;\;\sum_{c=0}^{n-(1-p)r} K_{2n+p,2\ell+q,2c+r} x^n y^\ell z^c
\end{align*}

Note that 
\begin{eqnarray*}
 S (x,y,z) &=& \sum_{p=0}^1 \sum_{q=0}^1 \sum_{r=0}^1 x^p y^q z^r \; S _{p,q,r}(x^2,y^2,z^2)
\end{eqnarray*}

After significant calculations, the resultant generating functions are:

\begin{align}
\label{GEN3P} x ( y - 2 ) (1 + z) S(x^2, y^2, z^2) = 
 2 (1  -  S(x, y, z)) + x y (1+z) S(x, y, z)^2 \\
\label{GEN3} S(x,y,z)={x (1+z) \over 2}\Biggl[  (2-y)\; S(x^2,y^2,z^2) + y\; S(x,y,z)^2\Biggr]+1
\end{align}

% cite Formal Power Series Paper. on ipad
As before if we formally assign the indeterminant value of $0^0$ to $1$ (this can also be handled with limits); then in that case we then have the generating function equation \ref{GEN3} yielding the known generating functions from equations \ref{GFE}, \ref{GEN1} and \ref{GEN2}.

When $z=0$ equation \ref{GEN3} yields equation \ref{GEN1}.

\begin{align*}
S(x,y,0)={x\over 2}\Biggl[  (2-y)\; S(x^2,y^2,0) + y\; S(x,y,0)^2\Biggr]+1
\end{align*}

When $y=1$ equation \ref{GEN3} yields equation \ref{GEN2}.

\begin{align*}
S(x,1,z)={x (1+z) \over 2}\Biggl[   S(x^2,1,z^2) + y\; S(x,1,z)^2\Biggr]+1
\end{align*}

When both $y=1$ and $z=0$ we have equation \ref{GEN3} yielding equation \ref{GFE}.

\begin{align*}
S(x,1,0)={x \over 2}\Biggl[   S(x^2,1,0) + y\; S(x,1,0)^2\Biggr]+1
\end{align*}

\nocite{riorcomb,RootIsoTree}

\bibliographystyle{elsarticle-harv}
\bibliography{main}




%% \bibitem must have one of the following forms:
%%   \bibitem[Jones et al.(1990)]{key}...
%%   \bibitem[Jones et al.(1990)Jones, Baker, and Williams]{key}...
%%   \bibitem[Jones et al., 1990]{key}...
%%   \bibitem[\protect\citeauthoryear{Jones, Baker, and Williams}{Jones
%%       et al.}{1990}]{key}...
%%   \bibitem[\protect\citeauthoryear{Jones et al.}{1990}]{key}...
%%   \bibitem[\protect\astroncite{Jones et al.}{1990}]{key}...
%%   \bibitem[\protect\citename{Jones et al., }1990]{key}...
%%   \harvarditem[Jones et al.]{Jones, Baker, and Williams}{1990}{key}...
%%


% \end{thebibliography}

\end{document}